\documentclass[12pt]{article}
\usepackage{amssymb}

\author{I. G. Korepanov\\
\normalsize Southern Ural State University\\[-0.5ex]
\normalsize 76 Lenin avenue\\[-0.5ex]
\normalsize 454080 Chelyabinsk, Russia\\[-0.5ex]
\normalsize E-mail: kig@susu.ac.ru}
\date{}
\title{Euclidean tetrahedra and knot invariants}

\def\be{\begin{equation}}
\def\ee{\end{equation}}

\def\pa#1#2{{\partial#1\over\partial#2}}

\def\Ker{\mathop{\rm Ker}\nolimits}
\def\myIm{\mathop{\rm Im}\nolimits}
\def\minor{\mathop{\rm minor}\nolimits}

\newtheorem{theorem}{Theorem}

\begin{document}
\maketitle

\begin{abstract}
We construct knot invariants on the basis of ascribing
Euclidean geometric values to a triangulation of sphere~$S^3$ where the
knot lies. The main new feature of this construction compared to the author's
earlier papers on manifold invariants is that now nonzero ``deficit
angles'' (in the terminology of Regge calculus) can also be handled. Moreover,
the knot goes exactly along those edges of triangulations that have nonzero
deficit angles.
\end{abstract}


\section{Introduction}

In this paper, we construct invariants of knots in the sphere~$S^3$ based,
very schematically, on the following ideas. A triangulation of~$S^3$ is
taken such that the knot goes along some of its edges. Then, Euclidean
geometric values are ascribed in some way to the elements of triangulation.
In particular, every edge gets a Euclidean lengths. Hence, the dihedral angles
also acquire some values, as well as (say) the volumes of tetrahedra. These
(dihedral angles and volumes) are taken with a minus sign for some tetrahedra,
according to certain orientation considerations.

Then, an acyclic complex is constructed in the spirit of
papers~\cite{Korepanov 2-4,Korepanov 1-5}, whose vector spaces consist
of the differentials of our geometric values. Recall that in those papers invariants of three-
and four-dimensional {\em manifolds\/} have been proposed based on algebraic
relations corresponding in a natural way to simplicial (Pachner) moves
--- elementary rebuildings of a manifold triangulation. The principal
new point in the present paper as compared with \cite{Korepanov 2-4,Korepanov 1-5}
(as well as~\cite{Korepanov JNMP,KM JNMP,Korepanov 3-3}) is that here an
algebraic relation is proposed dealing with the situation where the {\em
deficit angles\/} around some edges are nonzero modulo~$2\pi$. These deficit
angles are defined like in Regge calculus: minus sum of the dihedral angles
in all tetrahedra to which the edge belongs, save that now we take the {\em
algebraic\/} sum, with the signs mentioned in the previous paragraph. Construction
of a knot invariant follows as a natural application of such relation:
the edges in the triangulation along which the knot goes are singled out
exactly by the fact that the deficit angles for them are not zero but some
value~$(-\varphi)$.

Our knot invariant is expressed through the torsion of acyclic complex,
edge lengths and tetrahedron volumes much like the three-manifold 
invariant~\cite[formula~(5)]{Korepanov 2-4} (edges belonging to the knot
make, however, some difference, as we are going to explain).

Below, in section~\ref{sec moves} we present simplicial moves which are
enough to pass from a given triangulation of sphere~$S^3$ {\em with a knot
in it\/} to any other triangulation. In section~\ref{sec complex}, we present
the relevant acyclic complex and the formula for knot invariant. In section~\ref{sec 1to2}, we prove our new
key algebraic relation dealing with nonzero deficit angles. In the concluding section~\ref{sec discussion}
we discuss our results and prospects for further work.

\section{Simplicial moves}
\label{sec moves}

We are going to consider, strictly speaking, {\em pseudo}triangulations
of a sphere with a knot. From triangulations in the proper sense of the
word, they differ in that the boundary of a simplex in a pseudotriangulation
can contain {\em several times\/} the same simplex of a smaller dimension;
besides, a simplex in a pseudotriangulation is, generally, not determined
by the list of its vertices. Note that although it is often more convenient
for combinatorial topologists to consider triangulations in the proper sense 
(see especially the very useful paper~\cite{Lickorish}), it is
usually not very hard to bridge the way to pseudotriangulations.

Let there be a knot~$K$ in the sphere~$S^3$. Consider a pseudotriangulation
of this sphere obeying the following conditions:
\begin{itemize}
\item[(a)] the whole knot $K$ lies on some edges of the pseudotriangulation;
\item[(b)] for any tetrahedron in the pseudotriangulation, not more than
two of its vertices belong to~$K$;
\item[(c)] any edge~$a$ of the pseudotriangulation either has two {\em
different\/} vertices as its ends or, if its ends coincide, $a$ represents
a generator of the abelianization of the knot group~$\pi(K)$.
\end{itemize}

In other words, condition~(c) states that the edges with coinciding ends wind just
one time around the knot. We use such edges to describe ``moves $1\leftrightarrow
2$'' in the following theorem.

\begin{theorem}
\label{th1}
A pseudotriangulation of $S^3$ obeying the properties (a)--(c) can be transformed
into any other pseudotriangulation with the same properties by a sequence
of the following elementary moves:
\begin{itemize}
\item Pachner moves $2\leftrightarrow 3$ and $1\leftrightarrow 4$ (here,
for instance, ``$\,2\leftrightarrow 3$'' means that 2 tetrahedra are replaced
with 3 tetrahedra or back). Such moves must not affect the edges lying
on the knot~$K$ (the knot may, however, pass through edges and/or vertices
lying in the boundary of the transformed cluster of tetrahedra);
\item moves $1\leftrightarrow 2$ ``on the knot''. Let an edge~$BD$ lie
on the knot, and let there be a tetrahedron $BDAA$ in the pseudotriangulation,
with its edge $AA$ winding one time around the knot (cf.~the remark before
this theorem). The move $1\to 2$ is defined the following way: take a point~$C$
in the edge~$BD$ and replace the tetrahedron~$BDAA$ by two tetrahedra $BCAA$
and~$CDAA$. The move $2\to 1$ is the inverse to that.
\end{itemize}
\end{theorem}

{\it Proof\/} of this theorem is not very difficult for the careful reader
of paper~\cite{Lickorish}. It will be presented in a further publication.

\section{Acyclic complex}
\label{sec complex}

To a representation $f\colon \; \pi(K)\to G$ of the knot~$K$
group into some group~$G$, a branched over~$K$ covering of sphere~$S^3$
(where $K$ lies) naturally corresponds. Namely, take first  the {\em universal\/}
branched over~$K$ covering of~$S^3$ (where $\pi(K)$ naturally acts), and then identify its points $x$ and
$x'$ if and only if $x'$ is obtained from $x$ by the action of an element
of the subgroup $\Ker f \subset \pi(K)$. We will assume that $f$ is {\em nontrivial\/}:
$\myIm f\ne \{e\}$, where $e$ is the unit element of~$\pi(K)$.
Note that for $f$ to be nontrivial it is sufficient that the image of some
{\em overpass\/} of some knot~$K$ diagram, considered as an element of~$\pi(K)$,
should not equal~$e$ (then this will also hold for any overpass and any
diagram).

Below we confine ourselves to the case $G={\rm SO}(3)$. Then, the image
of any overpass is a rotation through (one and the same) angle~$\varphi$ about some axis
in the three-dimensional Euclidean space~$\mathbb R^3$, all these axes
going through the origin of coordinates. If there is only one such axis,
then $f$ is a {\em scalar\/} representation, i.e., a representation of
the abelianization of~$\pi(K)$.

We lift up the triangulation (or pseudotriangulation, we omit the prefix
`pseudo-' below) of~$S^3$ satisfying the properties (a)--(c)
 from section~\ref{sec moves} into the covering corresponding to representation~$f$.
Then we follow the same way as when constructing manifold invariants in
papers~\cite{KM JNMP,Korepanov 3-3,Korepanov 2-4,Korepanov 1-5}: to each
vertex of the lifted triangulation we ascribe coordinates in~$\mathbb R^3$
obeying the following condition: if vertices $F^{(1)}$ and~$F^{(2)}$ lie above
the same vertex~$F$ of the initial triangulation of~$S^3$ and $F^{(2)}=gF^{(1)}$,
where $g\in \pi(K)$, then the coordinates of $F^{(2)}$ are obtained from
those of~$F^{(1)}$ by the transformation $f(g)\in {\rm SO}(3)$. One more
condition is that of a ``general position'': any four vertices, if at least
two of them do not belong to the knot, must not
lie in the same (2-dimensional) plane (thus, the volumes of tetrahedra
considered below will be nonzero). Otherwise, the coordinates we ascribe
to vertices are arbitrary.

In such way, every edge of the {\em initial\/} triangulation acquires a Euclidean
length, and usual Euclidean values are also assigned to dihedral angles
and tetrahedron volumes. Some of these angles and volumes are taken, however,
with the minus sign according to the following rule.

We assume that all the tetrahedra of our triangulation of~$S^3$ are {\em oriented
consistently}, i.e., for every tetrahedron, an order of its vertices is
fixed up to even permutations; if there are two tetrahedra with a common
two-dimensional face $ABC$ and respective fourth vertices $D$ and~$E$,
then their consistent orientations will be, for instance, $ABCD$ and~$EABC$.
When we place in $\mathbb R^3$ an oriented tetrahedron $ABCD$, we take
as its volume the {\em oriented volume\/} $V_{ABCD}=\frac{1}{6}
\overrightarrow{AB} \overrightarrow{AC} \overrightarrow{AD}$ (triple scalar product).
Recall that $V_{ABCD}\ne 0$. So, the sign of the volume is already here,
and we assign to {\em all the dihedral angles of a tetrahedron\/} the same
sign as its volume has.

Here is our acyclic complex (explanations below):
\be
0\to \mathfrak a \stackrel{\left(\pa{x_i}{(z\;{\rm or}\;\varphi)}\right)}{\longrightarrow} (dx) 
\stackrel{\left(\pa{l_a}{x_i}\right)}{\longrightarrow} (dl) \stackrel{\left(\pa{\omega_b}{l_a}\right)}{\longrightarrow}
(d\omega) \to (\cdots) \to (\cdots) \to 0.
\label{complex}
\ee
It resembles very much the complex used in our construction of three-manifold 
invariants in~\cite[formula~(1)]{Korepanov 1-5}. Like in that paper, 
we are using somewhat loose but very convenient notations for vector
spaces consisting of differentials of Euclidean values.

Let us begin explanations from the space denoted $(dx)$. It consists of
column vectors whose entries are: three differentials $dx_i$, $dy_i$, $dz_i$
of the coordinates of every vertex~$i$ which does not lie on the knot, and one differential $dz'_k$ for every
vertex~$k$ that does lie on the knot. To be exact, we choose {\em one\/} preimage
in our lift-up for each vertex in the triangulation of~$S^3$, and denote its coordinates $x_i$,
$y_i$ and~$z_i$ if this $i$-th vertex does not lie on the knot, or take
just one coordinate~$z'_k$ if this $k$-th vertex does lie on the knot and
hence its preimage lies on some fixed axis~$z'$ in~$\mathbb R^3$ ($z'$ is of course one of
the axes mentioned in the second paragraph of this section).

The space $(dl)$ consists of column vectors $(dl_1,\ldots,dl_{N_1})^{\rm T}$,
where $l_a$ is the length of the $a$-th edge and $N_1$ is the number of
edges. It is clear how we define the mapping $(dx)\to (dl)$: it is the
{\em differential\/} of the mapping that sends the coordinates of vertices to
lengths of edges which join them.

Similarly, $(dl)\to (d\omega)$ is the differential of the mapping sending
the edge lengths to the deficit angles around the edges. By definition,
$\omega_b$ is minus algebraic sum of dihedral angles in all the tetrahedra
which contain the edge~$b$ (it is here that the angle signs come into play).
Here, the lengths $l_a$ are allowed to have arbitrary infinitesimal increments~$dl_a$,
so the values $\omega_b$ deviate by some $d\omega_b$ from either $0 \pmod{2\pi}$ or
$(-\varphi) \pmod{2\pi}$ for the respective edges lying or not lying on the knot.

Now we turn to the vector space $\mathfrak a$ and mapping $\mathfrak a\to
(dx)$. By definition, $\mathfrak a = \{0\}$ if $f$ is not a scalar representation.
In the case $f$ is scalar, we assume that $f$ sends $\pi(K)$ to rotations
about the axis~$z$. Then $\mathfrak a$ is, by definition, the two-dimensional
vector space consisting of columns denoted~$\pmatrix{dz\cr d\varphi}$.
The mapping $\mathfrak a\to (dx)$, by definition, adds $dz$ to the $z$-coordinates
of all vertices and rotates them all through the angle~$d\varphi$ around
axis~$z$.

In all these spaces $\mathfrak a$, $(dx)$, $(dl)$ and $(d\omega)$ there
are natural distinguished bases (up to the ordering of basis vectors).
Thus, all the mappings considered above are identified with their matrices.
Moreover, matrix $\left(\pa{\omega_b}{l_a}\right)$ is {\em 
symmetric}~\cite[section~2]{Korepanov JNMP}, so the remaining arrows in
the sequence~(\ref{complex}) are filled with matrices transposed to $\left(
\pa{l_a}{x_i}\right)$ and $\left( \pa{x_i}{(z{\rm\;or\;}\varphi)}\right)$,
with no notice to the geometric meaning of the spaces denoted $(\cdots)$,
in the same way as it has been done in~\cite{Korepanov 2-4,Korepanov 1-5}.

\begin{theorem}
The sequence (\ref{complex}) is an acyclic complex.
\label{th2}
\end{theorem}

{\it Proof}. The fact that the product of two consecutive arrows is zero,
follows immediately from geometric considerations. For instance, here is
how they go for the two arrows at the term~$(dx)$: if lengths $l_a$ are
calculated from the vertex coordinates~$x_i$ (and do not change arbitrarily
by themselves), then the deficit angles remain zero for edges not lying
on the knot and equal to $(-\varphi)$ --- for those lying there. On passing
to infinitesimals, we get at once that the product of these two arrows
is zero.

Thus we conclude that (\ref{complex}) is indeed an algebraic complex. As for the acyclicity,
it can be proved using methods of paper~\cite{Korepanov 2-4}, see section~2
of that work. The detailed proof will be presented in a further publication.

\medskip

Our knot invariant is defined by the formula
\be
I(K) = \tau \, \frac{\prod' l^2}{\prod 6V} \,
\Bigl(-2(1-\cos\varphi)\Bigr)^{N_0^{\rm knot}}.
\label{inv}
\ee
Here $\tau$ is the torsion of the complex~(\ref{complex}). In terms of
some minors of matrices (chosen according to the general theory of torsions),
and taking into account the symmetry of the complex, it can be written
as
\be
\tau = \frac{(\minor \left(\pa{l_a}{x_i}\right) )^2}{\minor \left(
\pa{\omega_b}{l_a}\right) (\minor \left(\pa{x_i}{z\;{\rm or}\;\varphi}\right)
)^2 };
\label{torsion}
\ee
the empty minor is considered to equal unity. The product in the
denominator in (\ref{inv}) is taken over all tetrahedra, while the primed product
in the numerator --- only over the edges {\em not lying on the knot}.
The number $N_0^{\rm knot}$ means the number of vertices in the triangulation
{\em lying on the knot}.

\begin{theorem}
The quantity I(K) is a knot invariant.
\label{th3}
\end{theorem}

{\it Proof\/} is based on investigating how the torsion behaves under the
simplicial moves from Theorem~\ref{th1}. For moves $2\leftrightarrow 3$
and $1\leftrightarrow 4$ this has been already done, essentially, 
in~\cite{Korepanov JNMP}, see also the parts of 
\cite{Korepanov 2-4,Korepanov 1-5} dealing with three-manifolds. As for
moves $1\leftrightarrow 2$, we are concerned with them in the next section.

\section{Algebraic relation for move $1\to 2$}
\label{sec 1to2}

Let a tetrahedron $BDAA$ be as in Theorem~\ref{th1}: edge~$BD$ lies on
the knot~$K$, while edge~$AA$ goes one time around~$K$. Let us have a representation
$f\colon \; \pi(K)\to {\rm SO}(3)$ such that the passing along $AA$ is
taken into the rotation about the axis~$z$ through an angle~$\varphi
\ne 0 \pmod{2\pi}$. When lifted to the covering corresponding to~$f$ (see
the first paragraph of section~\ref{sec complex}), the beginning and end
of edge~$AA$ become different points, denoted below as $A^{(1)}$ and~$A^{(2)}$.

Our aim now is to investigate how the torsion of complex~(\ref{complex})
behaves under the move $1\to 2$ described in Theorem~\ref{th1}. As the
form of our acyclic complex~(\ref{complex}) suggests, we should be interested in
vertices, edges and tetrahedra taken away or added to the simplicial
complex. Namely, the move $1\to 2$:
\begin{itemize}
\item adds new vertex $C$;
\item takes away the edge $BD$ and adds new edges $BC$, $CD$ and $CA$;
\item takes away tetrahedron $BDAA$ and adds new ones, $BCAA$ and $CDAA$.
\end{itemize}

Considerations related to the triangular form of matrices and similar to
those used in \cite{Korepanov JNMP} show that the factor by which the torsion of complex~(\ref{complex})
is multiplied under the move $1\to 2$ can be calculated by using the following ``local'' acyclic
complexes:
\be
0\to 0\to (dl_{BD}) \to (d\omega_{BD}) \to 0\to 0
\label{k1}
\ee
and
\begin{eqnarray}
&& 0\to (dz_C) \to (dl_{BC}, dl_{CD}, dl_{AC})\to (d\omega_{BC}, d\omega_{CD}, d\omega_{AC})
\nonumber\\
&& \quad \to (\hbox{the term symmetric to }dz_C) \to 0.
\label{k2}
\end{eqnarray}
Here, for instance, the term $(dz_C)$ means the one-dimensional space ---
just the differential of coordinate~$z_C$; $(dl_{BC}, dl_{CD}, dl_{AC})$
is a three-dimensional space, etc.

We denote the torsion of complex~(\ref{k1}) as $\tau_1$, and the torsion of complex~(\ref{k2})
--- as~$\tau_2$. Direct calculation shows that
\be
\frac{\tau_2}{\tau_1} = \frac{6V_{CAAB}\cdot 6V_{DAAC}}{2 (1-\cos \varphi )\cdot
l_{AC}^2 \cdot 6V_{DAAB}}
\label{k3}
\ee
(and (\ref{k3}) also gives the ratio of torsions of the complexes corresponding
to the whole triangulations of~$S^3$ differing in the $1\to 2$ move). The
meaning of the result~(\ref{k3}) is the following: the torsion~(\ref{torsion})
of complex~(\ref{complex}):
\begin{itemize}
\item acquires the factor $-\frac{1}{2(1-\cos \varphi )}$ on adding a vertex ($C$ in
our case) lying on the knot;
\item acquires the factor $l_i^{-2}$ on adding an edge~$i$ not lying
on the knot. Note that for edges {\em lying\/} on the knot, there are no
such factors;
\item acquires the factor $6V$ for every appearing tetrahedron (of
volume~$V$), and the factor $(6V)^{-1}$ for every disappearing tetrahedron.
\end{itemize}
This all together guarantees the invariance of value~$I(K)$ under moves
$1\leftrightarrow 2$, which concludes the proof of Theorem~\ref{th3}.

\section{Discussion of results}
\label{sec discussion}

The main thing in this paper is formula~(\ref{k3}). It shows that a relation
corresponding to a relevant simplicial move does exist for a nonzero deficit
angle as well. In the present paper, this fact is applied to knot theory.
On the other hand, it may find applications also in Regge calculus, i.e.,
discrete gravity theory.

It may make sense to search for a quantum analogue of formula~(\ref{k3}).
This problem seems to be solvable in our case of invariants built using
Euclidean geometry in 3-simplices, because of analogy with the relation
for move $2\to 3$, which played the key role in constructing {\em manifold\/}
invariants, started in~\cite{Korepanov JNMP}, and which can be obtained
in a (double) semiclassical limit from the pentagon equation for quantum 6$j$-symbols.

One more apparently solvable problem is generalization of our geometric
constructions onto the case of multidimensional analogues of knots, and
onto the case where we apply some other than Euclidean geometry. Much more
difficult looks the problem of {\em quantizing\/} of these latter constructions.

The results of actual calculations of our invariants $I(K)$
for specific knots~$K$ and representations of group~$\pi(K)$ will be presented
in E.V.~Martyushev's work~\cite{Martyushev}. Also, in a future paper of
larger size, complete proofs of our theorems will be given. Here we just mention that
$$
I(\hbox{unknot})= - 4(1-\cos\varphi)^2,
$$
(which is not hard to calculate even without a computer) and that calculations
show a connection with Alexander polynomial and  Reidemeister torsion (both
abelian and nonabelian). Most interesting results are expected for most
complicated representations of~$\pi(K)$ in some group of geometric origin, as the analogy with our manifold
invariants~\cite{Martyushev manifolds 1,Martyushev manifolds 2} suggests.

\bigskip

{\bf Acknowledgements. }\vadjust{\nobreak}The work has been performed with a partial financial
support from Russian Foundation for Basic Research under Grant no.~04-01-96010.
I thank Evgeniy Martyushev for the careful reading of this paper and corrections.

\end{document}